\documentstyle[12pt]{article}
\title{ON THE AVERAGE NUMBER OF SHARP  CROSSINGS OF  CERTAIN GAUSSIAN   RANDOM  POLYNOMIALS}
\author{S. Shemehsavar,  S. Rezakhah\thanks{ \scriptsize Faculty of Mathematics and Computer Science, Amirkabir University of Technology, 424
 Hafez Avenue, Tehran Iran. Email: Rezakhah@aut.ac.ir,$\;\;$ Shemehsavar@aut.ac.ir.}}
\date{}

\begin{document}

\maketitle

\pagestyle{plain}

\renewcommand{\theequation}{\arabic{section}.\arabic{equation}}

\begin{abstract}
 Let  $Q_n(x)=\sum _{i=0}^{n} A_{i}x^{i}$ be a random algebraic polynomial where
 the coefficients
  $A_0,A_1,\cdots $ form a sequence of centered Gaussian random variables. Moreover,
  assume that
the increments $\Delta_j=A_j-A_{j-1}$, $j=0,1,2,\cdots$ are
independent, assuming $A_{-1}=0$. The coefficients can be
considered as $n$ consecutive observations of a Brownian motion.
We obtain the asymptotic behaviour of the expected number of
u-sharp crossings
 of polynomial $Q_n(x)$ . We refer to u-sharp
crossings as  those zero up-crossings with slope greater than $u$,
or those
 down-crossings  with slope smaller than $-u$.  We consider the
 cases
 where $u$ is unbounded and is  increasing with $n$, where
 $u=o(n^{5/4})$, and $u=o(n^{3/2})$ separately.

\end{abstract}

\noindent Keywords and Phrases: random algebraic polynomial,
number of real zeros, sharp crossings, expected density, Brownian   motion.\\
AMS(2000) subject classifications. Primary 60H42, Secondary 60G99.

 \setcounter{equation}{0}
\newcommand{\be}{\begin{equation}} \newcommand{\ee}{\end{equation}}
\section{Preliminaries}
\hspace{.2in}
 The theory of the expected number of real zeros of
random algebraic polynomials was addressed in  the fundamental
work of M. Kac\cite{kac1}.
 The works of Logan and Shepp \cite{lsh1,lsh2}, Ibragimov and
 Maslova \cite{ibm}, Wilkins \cite{wil}, and Farahmand
 \cite{far1} and Sambandham \cite{sam1,sam2} are other fundamental
 contributions  to the subject. For various aspects on random
 polynomials see Bharucha-Reid and Sambandham \cite{bhs}, and
 Farahmand \cite{far2}.
    There has been recent interest in cases where the
coefficients form certain random processes,
  Rezakhah and Soltani \cite{rs1,rs2}, Rezakhah and Shemehsavar \cite{rs3}.

 Let
$$ Q_{n}(x)=\sum_{i=0}^{n}A_{i}x^{i},\, \; -\infty<x<\infty , \eqno{(1.1)}$$
where  $A_{0},A_{1},\cdots  $ are mean zero Gaussian random
variables for which the increments  $\Delta _i =A_i -A_{i-1}, \;
i=1,2,\cdots $ are independent, $A_{-1}=0$. The sequence $A_0,A_1
\cdots $
 may be considered as successive Brownian points,
 i.e., $A_j=W(t_j),\; j=0,1,\cdots $, where $t_0<t_1< \cdots $
and $\{W(t),\;t\geq 0\}$ is the standard Brownian motion.
In this physical interpretation, Var$(\Delta_j)$ is the distance
 between successive times $t_{j-1},\; t_j$.
    Thus  for $j=1,2,\cdots$, we have that $A_j=\Delta_0 +\Delta _1
+\cdots +\Delta_j$, where $\Delta_i \sim N(0,\sigma_i^2 )$ and
are independent.
  Thus $
Q_n(x)=\sum_{k=0}^n (\sum _{j=k}^n x^j)\Delta_k=\sum_{k=0}^n a_k(x)\Delta _k,
$
and $
  Q'_n(x)=\sum_{k=0}^n (\sum _{j=k}^n jx^{j-1})\Delta_k=\sum_{k=0}^n b_k (x)\Delta _k,
$
where
    $$a_k (x) = \sum _{j=k}^n x^j,\;\;
 \; \; \; b_k (x) = \sum _{j=k}^n
jx^{j-1},\;\;\hspace{.3in}  k=0,\cdots , n.\eqno{(1.2)}$$

We say  that   $Q_n(x)$ has a zero up-crossing at $t_{0}$ if
there exists $\varepsilon>0$ such that $Q_n(x)\leq 0$ in
$({t_{0}-\varepsilon,t_{0}})$ and $Q_n(x)>0$ in
$(t_{0},t_{0}+\varepsilon)$. Similarly $Q_n(x)$ is said to have
 a zero down crossing at $t_0$, if
  $Q_n(x)>0$ in
$({t_{0}-\varepsilon,t_{0}})$ and $Q_n(x)\leq 0$
  in
$(t_{0},t_{0}+\varepsilon)$. We  study the asymptotic behavior of
the expected number of, $u$-sharp zero crossings,  those zero
up-crossings with slope greater than $u>0$, or  those
down-crossings with slope smaller than $-u$.

 Cramer and Leadbetter
[1967 p 287]  have shown that the expected number of total zeros
of any Gaussian non stationary process, say  $Q_n(x)$,  is
calculated via the following formula
$$
EN(a,b)=\int_a^b dt \int_{-\infty}^{\infty}
  |y|p_{t}(0,y)dy
$$
where $p_{t}(z,y)$ denotes the joint  density of  $Q_n(x)=z$ and
its derivative  $Q'_n(x)=y$, and
\begin{eqnarray*}
p_{t}(0,y)\!\!&\!=\!&\!\![2\pi\gamma \sigma (1-\mu^2)^1/2]^{-1}
\\
&&\times \exp\bigg \{-\frac{\gamma^2m^2+2\mu\gamma \sigma
m(y-m'+\sigma^2(y-m')^2}{2\gamma^2 \sigma^2(1-\mu^2)} \bigg \}
\end{eqnarray*}
in which $m=E(Q_n(x))$, $\;m'=E(Q'_n(x))$, $\;
\sigma^2=\mbox{Var}(Q_n(x))$, $\;\gamma^2=\mbox{Var}(Q'_n(x))$,
and $\mu=\mbox{Cov}(Q_n(x),Q'_n(x))/(\gamma \sigma)$.

By a similar method as in Farahmand \cite{far2} we  find  that
$ES_u(a,b)$, the expected number of $u$-sharp zero crossings   of
$Q_n(x)$  in any interval $(a,b)$, satisfies
$$ES_u(a,b)=\int_a^b dt \left\{\int_{-\infty}^{-u} +\int_u^{\infty}\right \}
  |y|p_{t}(0,y)dy=
\int_a^b f_n(x)dx\eqno{(1.3)}$$  where
$$ f_n(x)=\frac{1}{\pi}g_{1,n}(x)\exp(g_{2,n}(x))\eqno{(1.4)}$$
and
$$ g_{1,n}(x)= E A^{-2},\;\;\;\;\;\;\;\;\;\;\;\;
g_{2,n}(x)=-\frac{A^2u^2}{2E^2},\;\;
 \eqno{(1.5)}$$
in which
$$\;
A^2=\mbox{Var}(Q_n(x))\;=\;\sum_{k=0}^na_k^2(x)\sigma^2_k,\;\;\;\;
B^2=\mbox{Var}(Q'_n(x))=\sum_{k=0}^nb_k^2(x)\sigma^2_k,\;$$
$$
D\!=\!\mbox{Cov}(Q_n(x), Q'_n(x))\!=\!\sum_{k=0}^n
a_k(x)b_k(x)\sigma^2_k,  \;\;\;\; \mbox{and} \;\;\;\;
E^2=A^2B^2-D^2,$$ and $a_k(x),\;$  $\;b_k(x)$ is defined by (1.2).

\section{Asymptotic behaviour of $ES_u$ }
In this section we obtain  the asymptotic behaviour of the
expected number of
 $u$-sharp crossings of  $Q_{n}(x)=0$  given by (1.1).
 We prove  the following theorem for the case that
 the increments
 $\Delta_1 \cdots \Delta_n$ are independent and  have the same distribution. Also we assume that
  $\sigma_k^2=1$, for $k=1\cdots n$.\\

  \noindent
 \textbf{Theoream(2.1)}: Let $Q_n(x)$ be the random algebric
 polynomial given by (1.1) for which
 $A_{j}=\Delta_{1}+...+\Delta_{j}$ where $\Delta_{k}$, $k=1,...,n$
 are independent and $\Delta_{j}\sim N(0,\sigma^2_{j})$ then the
 expected number of u-sharp crossings of $Q_n(x)=0$  satisfies

\noindent
 (i)- for $u=o(n^{5/4})$
\begin{eqnarray*}
ES_{u}(-\infty,\infty)\!&\!=\!&\!\frac{1}{\pi}\log(2n+1)+\frac{1}{\pi}(1.920134502)\\&&+\frac{1}{\pi
\sqrt{2n}}\left( - \pi +2 \arctan ( \frac {1}{2\sqrt{2n}})\right)
+\frac{C1}{n\pi}\\&&+\frac{u^2}{n^3\pi}
\left(19.05803659-\frac{8}{3}\ln(n^3+1)\right)\\&&
+\frac{u^4}{n^6\pi}\left(-34989.96324+\frac{32}{3}\ln(n^6+1)\right)+o(n^{-1})
\end{eqnarray*}
where $C1=-0.7190843756$ for $n$ even and $C1=1.716159410$ for $n$
odd.

\noindent (ii)-
 for $u=o(n^{3/2})$
$$ES_{u}(-\infty,\infty)=\frac{1}{\pi}\log(2n+1)+\frac{1}{\pi}(1.920134502)+o(1)$$

\noindent
 \textbf{proof:} Due to the behaviour of $Q_n(x)$, the
asymptotic behaviour is treated separately on
  the intervals $1<x<\infty$,
$-\infty<x<-1$, $0<x<1$ and $-1<x<0$.    For $1<x<\infty$, using
(1.3),  the change of variable $x=1+\frac{t}{n}$ and the equality
  $\left( 1+\frac{t}{n}
\right)^n=e^t\left(1-\frac{t^2}{n} \right) +O\left(\frac{1}{n^2}
\right)$, we find  that $$ES_u(1 ,\infty )=\frac{1}{n}
\int_{0}^{\infty}f_n\big(1+\frac{t}{n}\big)dt,\;\;\;
$$ where  $f_n(\cdot )$ is defined by (1.4). Using  (1.5), and by tedious manipulation we have that
$$g_{2,n}\left(1+\frac{t}{n}\right)=o(n^{-2})$$
 and
$$
n^{-1}g_{1,n}\left(1+\frac{t}{n}\right)=
\left(R_1(t)+\frac{S_1(t)}{n}+O\left(
\frac{1}{n^2}\right)\right),\;\;\;   \eqno(2.1)
$$
where {\small
$$
R_1(t)=\frac
{\sqrt{(4t\!-\!15)e^{4\,t}\!+\!(24t+32)e^{3\,t}\!-\!e^{2\,t}(8t^
3\!+\!12t^2\!+\!36t\!+\!18) +8e^{t}t+1}}{2\,t\left
(-1-3\,{e^{2\,t}}+4\,{e^ {t}}+2\,t{e^{2\,t}}\right )}
$$
and  $S_1(t)=S_{11}(t)/S_{12}(t)$ in which
\begin{eqnarray*}
S_{11}(t)\!\!&\!\!=\!\!&\!\! -
0.25\left((4t^2-6t-27)e^{6t}+(156-84t+116t^2-24t^3)e^{5t}
\right.\hspace{1.4in}\\&&
 \!+( 16t^5\!- \!72t^4+ 96t^3- 212t^2\!+\! 220t-
331)e^{ 4t} \!+\! (328\!-\!168t+128t^2\!-\!104t^3)e^{3t} \\&&
\left. +( 8t^4 + 8t^3-32t^2+ 42t- 153 )e^{2t}+(28-4t-4t^2)e^t
-1\right)
\end{eqnarray*}
\begin{eqnarray*}
S_{12}(t)\!&\!=\!&\!\left(
2\,{e^{2\,t}}t-1+4\,{e^{t}}-3\,{e^{2\,t}} \right) ^{2
}\\&&\hspace{-.25in}\times({{1-8
{t}^{3}{e^{2t}}\!-\!12{e^{2t}}{t}^{2}+8{e^{t}}t\!-\!18{e^{2t}}\!-\!36
{e^{2t}}t\!-\!15{e^{4t}}+32{e^{3t}} +24{e^{3t}}t+4t{e^{4
t}}}})^{1/2}
\end{eqnarray*}
 One can easily verify that as
$\;t\rightarrow \infty$,
$$
R_1(t)=  \frac{1}{2t^{3/2}}+O(t^{-2}),\;\;\;\; \hspace{1in}
S_1(t)= - \frac{1}{8t^{1/2}}+O(t^{-3/2}).
$$
As (2.1) can not be integrated term by term,
  we use the equality
$$
\frac{I_{[t>1]}}{8n\sqrt{t}}=\frac{I_{[t>1]}}{8n\sqrt{t}+t\sqrt{t}}
+O\left(\frac{1}{n^2}\right),  \eqno(2.2)
$$where
$$I_{[t>1]}=\left\{\begin{array}{ll}1\;\;\;\; \mbox{if }\;\;t\geq 1\\0\;\;\;\;\mbox{if} \;\; t<1\end{array}\right. ,
$$
 Thus by (2.2) we have that
$$
\frac{1}{n}f_n\bigg(1+\frac{t}{n}\bigg)=-\frac{I_{[t>1]}}{\pi(8n\sqrt{t}+t\sqrt{t})}+\frac{R_1(t)}{\pi}
+\frac{1}{\pi} \left( \frac{S_1(t)}{n}+
\frac{I_{[t>1]}}{8n\sqrt{t}}\right)+O\left(\frac{1}{n^{2}}\right).
$$
This expression is term by term integrable, and provides that
\begin{eqnarray*}ES_u(1,\infty)\!&\!=\!&\! \frac{1}{n}\int_0^{\infty} f_n\bigg(1+\frac{t}{n}\bigg)dt  =
\frac{1}{2\pi\sqrt{2n}} \left( - \pi +2
\arctan ( \frac {1}{2\sqrt{2n}})\right)\\
&&\;\;\; +\frac{1}{\pi}\int_0^{\infty}R_1(t)dt+\frac{1}{\pi n}
\int_0^{\infty} \left( S_1(t)+
\frac{I_{[t>1]}}{8\sqrt{t}}\right)dt+O\left(\frac{1}{n^{2}}\right),
\end{eqnarray*}
 where $\;\; \int_0^{\infty}R_1(t)dt= 0.7348742023\;$,
  and $\;\;\int_0^{\infty} \left( S_1(t)+
  I_{[t>1]}(8\sqrt{t})^{-1}\right)dt=-0.2496371198$.

For  $-\infty<x<-1$, using  $x=-1-\frac{t}{n}$,
  we have
$ES_u(-\infty, -1 )=\frac{1}{n} \int_{0}^{\infty}
f_n\left(-1-\frac{t}{n}\right)dt,\,$  where $f_n(\cdot )$ is
defined by (1.4), and  by (1.5)
$$
n^{-1}g_{1,n}\left(-1-\frac{t}{n}
\right)=\left(R_2(t)+\frac{S_2(t)}{n}
+O\left(\frac{1}{n^2}\right)\right)\eqno(2.3)
$$
where {\small
$$R_2(t)=   1/2\sqrt {{\frac {-2\,{e^{2\,t}}+{e^{4\,t}}+1-12\,{e^{2\,t}}{t}^{2}-8\,{
t}^{3}{e^{2\,t}}+4\,t{e^{4\,t}}-4\,t{e^{2\,t}}}{{t}^{2}\left( {e^
{2\,t}}-1+2\,t{e^{2\,t}} \right) ^{2}}}};
$$
$$\hspace{-.07in}\mbox{ for $n$ even}\;\;\;\;\;\; S_{2}(t)=\frac{S_{21}(t)+S_{22}(t)}{4S_{23}(t)} ,
\;\;\;\;\;\;\mbox{and  for $n$
odd}\;\;\;\;\;\;\;S_{2}(t)=\frac{S_{21}(t)-S_{22}(t)}{4S_{23}(t)}$$
in which
\begin{eqnarray*}
S_{21}(t)\!\!&\!=\!\!&\!  1+\left(
-8\,{t}^{4}+30\,t-8\,{t}^{3}+48\,{t}^{2}-3
 \right) {e^{2\,t}} \\&& + \left( 3-12t+52t^2+96t^3+40t^{4}-16\,{t}^{
5} \right) {e^{4\,t}}\! - \left( 18t+4t^{2}+1 \right) {e^{6\,t}},
\\
S_{22}(t)\!\!&\!=\!\!&\!  \left( 8\,t+32\,{t}^{3}+40\,{t}^{2}
 \right) {e^{3\,t}}  +
 \left( -8\,{t}^{2}-12\,t \right) {e^{5\,t}}+4\,{e^{t}}t,
 \\ S_{23}(t)\!\!&\!=\!\!&\!  \left(e^{4\,t}(4t+1)-2e^{2t}(1+2t+6t^{2}+4t^{3})+1\right)^{1/2}
  \left( e^{2\,t}(2t+1)-1
 \right) ^{2}.
\end{eqnarray*}
and
$$
g_{2,n}\left(-1-\frac{t}{n}\right)=\left(\frac{u^2}{n^3}\right)G_{2,1}(t)+
o\left(\frac{1}{n}\right)\eqno(2.4)
$$
where
$$
G_{2,1}(t)={\frac { 16\left( 2\,{e^{2\,t}}t+{e^{2\,t}}-1 \right)
{t}^{3}}{2\,{e
^{2\,t}}-{e^{4\,t}}+12\,{e^{2\,t}}{t}^{2}-1+8\,{e^{2\,t}}{t}^{3}-4\,{e
^{4\,t}}t+4\,{e^{2\,t}}t}}
$$
 It  can be seen that as $t \rightarrow \infty$,
$$
 R_2(t)= \frac{1}{2t^{3/2}}+O(t^{-2}),\;\;\; \;\;\;
S_2(t)= \frac{-1}{8t^{1/2}}+O(t^{-3/2}), \;\;\;\;\;\;\;
G_{2,1}(t)=o(e^{-t}).$$
 Now by using the equality (2.2) we have that
\begin{eqnarray*}
\frac{1}{n}f_n(-1-\frac{t}{n})\! &\!
=\!&\!-\frac{I_{[t>1]}}{\pi(8n\sqrt{t}+t\sqrt{t})}
+\frac{1}{\pi}\left(R_2(t) +\frac{1}{n}(
S_2(t)+\frac{I_{[t>1]}}{8\sqrt{t}})\right.\\&&+\frac{u^2}{n^3}R_2(t)G_{2,1}(t)
+\frac{u^4}{2n^6}R_2(t)G^2_{2,1}(t)\bigg) +o\left(n^{-1}\right).
\end{eqnarray*}
Thus
\begin{eqnarray*}
ES_u(-\infty,-1)\!&\!=\!&\!
\frac{1}{n}\int_0^{\infty}f_n(-1-\frac{t}{n}) =\frac{1}{2\pi
\sqrt{2n}} \left( - \pi +2 \arctan ( \frac
{1}{2\sqrt{2n}})\right)\\&&+ \frac{1}{\pi}\int_0^{\infty}R_2(t) dt
+\frac{1}{n\pi}\int_0^{\infty}\left(
S_2(t)+\frac{I_{[t>1]}}{8\sqrt{t}}\right)dt\\&&
+\frac{u^2}{n^3\pi}\int_0^{\infty}R_2(t)G_{2,1}(t)dt
+\frac{u^4}{n^6\pi}\int_0^{\infty}\frac{1}{2}R_2(t)G^2_{2,1}(t)dt
+o\left(n^{-1}\right)
\end{eqnarray*}
where $\int_0^{\infty}R_2(t)dt=1.095640061$ ,
 and
$$ \int_0^{\infty}R_2(t)G_{2,1}(t)dt=-2.418589510,
\hspace{.2in}\int_0^{\infty}\frac{1}{2}R_2(t)G^2_{2,1}(t)dt=7.057233216,$$
 and for $n$ odd $\;\;\; \int_0^{\infty}\left(
S_2(t)+I_{[t>1]}(8\sqrt{t})^{-1}\right)dt =-0.0322863105 ,\;\;\;$
and for $n$ even  \\$\; \int_0^{\infty}\left( S_2(t)
+I_{[t>1]}(8\sqrt{t})^{-1}\right)dt =-0.4677136959.$

For $0<x<1$, let $x=1-\frac{t}{n+t}$, then $ES_u(0,1
)=\left(\frac{n}{(n+t)^2} \right) \int_{0}^{\infty} f_n
\left(1-\frac{t}{n+t}\right)dt,\;\;\;$ where $f_n(\cdot )$ is
defined by (1.4), and  by (1.5)
$$
g_{2,n}\left(1-\frac{t}{n+t}\right)=o(n^{-2}),\eqno(2.5)
$$
and
\begin{eqnarray*}
\left(\frac{n}{(n+t)^2}
\right)g_{1,n}\left(1-\frac{t}{n+t}\right)&\!=\!&\left(1-\frac{2t}{n}
+O\left( \frac{1}{n^2} \right)\right) \left(R_3(t)
+\frac{S_3(t)}{n}+O\left(\frac{1}{n^2}\right)\right)\\
& =\! & \left(R_3(t)+\frac{S_3(t)-2tR_3(t)}{n}\right)
+O\left(\frac{1}{n^2}\right), \;\;\;\;\;\;(2.6)
\end{eqnarray*}
where we observe that {\small $R_3(t) \equiv  R_1(-t)$
 and
 $S_3(t)=S_{31}(t)/S_{32}(t)$,  in which
 \begin{eqnarray*}
S_{31}(t)\!&\!=\!&\!\left(  \left( -7{t}^{2}-{\frac {69}{2}}t-{
\frac {63}{4}} \right) {e^{-6t}}+ \left( 6{t}^{3}+35t-55{t}^{2}+39
 \right) {e^{-5t}}\right.\\&&\!\!\!+ \!\! \left( 49t-4{t}^{5}+22{t}^
{4}+91{t}^{2}-{\frac {63}{4}}-12{t}^{3} \right) {e^{-4t}}+ \left(
-6{t}^{3}-30-44{t}^{2}-66t \right) {e^{-3t}}\\&&+\left.\left(
{\frac {35}{2}}t+2{t}^{4}-6{t}^{3}+16{t}^{2}+{\frac {123}{4}}
\right) {e^{-2t}}
 + \left( -
9-t-{t}^{2} \right) {e^{-t}}+3/4  \right),
\end{eqnarray*}
and $ S_{32}(t)\equiv S_{12}(-t)$}. Now as $t\to\infty\;,$  $
R_3(t)= \frac{1}{2t}+O(t^{-1/2}e^{-t/2}),\;$ and $ \; S_3(t)=
\frac{3}{4}+O(t^2e^{-t}) $.

 Since the relation (2.6) is not term by term integrable we use the
 equality
$$\frac{I_{[t>1]}}{2t}-\frac{I_{[t>1]}}{4n+2t}=\frac{I_{[t>1]}}{2t}-
\frac{I_{[t>1]}}{4n}+O\left(\frac{1}{n^2}\right),\;\eqno(2.7)$$
  can be written as
\begin{eqnarray*}
\frac{n}{(n+t)^2}f_n(1-\frac{t}{n+t})\! &\! =\!&\!
\frac{1}{\pi}\left(R_3(t)-\frac{I_{[t>1]}}{2t}\right)
+\frac{1}{\pi} \left(
 \frac{I_{[t>1]}}{2t}-\frac{I_{[t>1]}}{4n+2t}\right)
\\&&+\frac{1}{n\pi}\left( S_3(t)-2tR_3(t)+\frac{I_{[t>1]}}{4}\right)
+O\left(\frac{1}{n^2}\right)
\end{eqnarray*}
Thus we have that
\begin{eqnarray*}
ES_u(0,1)\!&\!=\!&\!\frac{n}{(n+t)^2}\int_0^{\infty}f_n(1-\frac{t}{n+t})dt\\&&=
\frac{1}{\pi}\int_0^{\infty}\left(R_3(t)-\frac{I_{[t>1]}}{2t}\right)dt+\frac{1}{2\pi}
(
 \log (4n+2)-\log (2))
\\&&+\frac{1}{n\pi}\int_0^{\infty}\left( S_3(t)-2tR_3(t)+\frac{I_{[t>1]}}{4}\right)dt
 +O\left(\frac{1}{n^2}\right) \;\;\;\;\;(2.8)
\end{eqnarray*}
where
$\int_0^{\infty}\left(R_3(t)-\frac{I_{[t>1]}}{2t}\right)dt=-0.2897712456\;\;\;$
and $\;\;\;\int_0^{\infty}\left(
S_3(t)-2tR_3(t)+\frac{I_{[t>1]}}{4}\right)dt=0.498174649.$

For $-1<x<0$, let $x=-1+\frac{t}{n+t}$, then $ES_u(-1,0 )=\left(
\frac{n}{(n+t)^2}\right ) \int_{0}^{\infty}
f_n\left(-1+\frac{t}{n+t}\right)dt,\;$ again by (1.4) and (1.5) we
have
$$
\left( \frac{n}{(n+t)^2}\right )g_{1,n}(t)=\left(
\frac{n^2}{(n+t)^2}\right)
\left(R_4(t)+\frac{S_4(t)}{n}+O\left(\frac{1}{n^2}\right)
\right),\eqno(2.9)
$$
in which {\small $R_4(t)\equiv R_2(-t)$.
$$\hspace{-.07in}\mbox{ for $n$ even}\;\;\;\;\;\; S_{4}(t)=\frac{S_{41}(t)+S_{42}(t)}{4S_{43}(t)} ,
\;\;\;\;\;\;\mbox{and  for $n$
odd}\;\;\;\;\;\;\;S_{4}(t)=\frac{S_{41}(t)-S_{42}(t)}{4S_{43}(t)}$$
 In which
 \begin{eqnarray*}
S_{41}(t)\!\! &\!=\!&\!\! 8\, \left(
-9/4\,t+6\,{t}^{2}-3\,{t}^{3}-{\frac {9}{8}}+{t}^{4} \right)
{e^{-2\,t}}\\&& +8\, \left( 15\,{t}^{4}-3/2\,t-22\,{t}^{3}+{ \frac
{9}{8}}+19/2\,{t}^{2}-2\,{t}^{5} \right) {e^{-4\,t}}\\&& +8\,
\left( {\frac {15}{4}}\,t-7/2\,{t}^{2}-3/8 \right) {e^{-6\,t}}+3,
\end{eqnarray*}
 $S_{42}(t)\equiv S_{22}(-t)$,   and
$ S_{43}(t)\equiv S_{23}(-t).\;$
Finally we have
$$g_{2,n}\left(-1+\frac{t}{n+t} \right)=
 \left(\frac{u^2}{n^3}\right) G_{2,1}(t)+o(n^{-1})\eqno(2.10)
$$
where
$$
G_{2,1}(t)=-{\frac {16 \left( -{e^{-2\,t}}+1+2\,{e^{-2\,t}}t
\right) {t}^{3}}{(-
4\,{e^{-4\,t}}t+{e^{-4\,t}}-12\,{e^{-2\,t}}{t}^{2}+4\,{e^{-2\,t}}t+1+8
\,{e^{-2\,t}}{t}^{3}-2\,{e^{-2\,t}})}}
$$

 As $t\rightarrow \infty $ we
have
$$R_4(t)=\frac{1}{2t}+O(t^{1/2}e^{-t}) ,\;\;\;\;\;\;
S_{4}(t)=\frac{3}{4}+O(te^{-t}),\;\;\;\;\;\;
G_{2,1}(t)=-16t^3+O(t^4e^{-2t}).
$$
We have that
\begin{eqnarray*}
\frac{n^2}{(n+t)^2} f_n(-1+\frac{t}{n+t})&\!=\!&\frac{1}{\pi}
\left(1-\frac{2t}{n} +O(\frac{1}{n^2})\right)
\left(R_4(t)+\frac{S_4(t)}{n}+O(\frac{1}{n^2})\right)\\&&
\times\left(1+\frac{u^2}{n^3}G_{2,1}(t)+\frac{u^4}{2n^6}G^2_{2,1}(t)
+o(n^{-1})\right)\\&&=
\frac{1}{\pi}\bigg\{R_4(t)+\frac{S_4(t)-2tR_4(t)}{n}+\frac{u^2}{n^3}R_4(t)G_{2,1}(t)\\&&
\hspace{.8in}+\frac{u^4}{2n^6}R_4(t)G^2_{2,1}(t) \bigg\}+o(n^{-1})
\end{eqnarray*}
Since this is not term by term integrable we use (2.7) and the
following equalities to solve this
$$\frac{8u^{2}t^{2}}{n^3}=\frac{8u^2t^2}{n^3+exp(t^3)}+o(n^{-2})\; ,\;\;\;\;\;\;
\frac{64u^4t^5}{n^6}=\frac{64u^4t^5}{n^6+exp(t^6)}+o(n^{-2})$$
Thus we have that
\begin{eqnarray*}
ES_u(-1,0)\!&\!=\!&\!\frac{n^2}{(n+t)^2}\int_0^{\infty}f_n(-1+\frac{t}{n+t})\\&&=\frac{1}{2\pi}
(
 \log (2n+1))+\frac{1}{\pi}
\left\{\int_0^{\infty}(R_4(t)-\frac{I_{[t>1]}}{2t})dt
(2))dt\right.\\&&+\frac{1}{n}\int_0^{\infty}(S_4(t)-2tR_4(t)+\frac{I_{[t>1]}}{4})dt\\&&
+\frac{u^2}{n^3}\int_0^{\infty}(R_4(t)G_{2,1}(t)+8t^2)dt
-u^2\int_0^{\infty}\frac{8t^2}{n^3+exp(t^3)}dt\\&&
+\frac{u^4}{n^6}\int_0^{\infty}(\frac{1}{2}R_4(t)G^2_{2,1}(t)-64t^5)dt
+u^4\int_0^{\infty}\frac{64t^5}{n^6+exp(t^6)}dt\bigg\}\\&&
+o(n^{-1})\;\;\;\;\;\hspace{1.5in}(2.11)
\end{eqnarray*}
where$\int_0^{\infty}(R_4(t)-\frac{I_{[t>1]}}{2t})dt=.3793914851,$
and$$\int_0^{\infty}(R_4(t)G_{2,1}(t)+8t^2)dt=21.47662610,\hspace{.2in}
\int_0^{\infty}(\frac{1}{2}R_4(t)G^2_{2,1}(t)-64t^5)dt=-34997.02047
 $$
 and we have for $n$ even
$\int_0^{\infty}(S_4(t)-2tR_4(t)+\frac{I_{[t>1]}}{4})dt=-0.4999081999,\;$
and for $n$ odd
$\int_0^{\infty}(S_4(t)-2tR_4(t)+\frac{I_{[t>1]}}{4})dt=1.499908200.
$ Also we have  that
$$\int_0^{\infty}\frac{8t^2}{n^3+exp(t^3)}dt=\frac{8ln(n^3+1)}{3n^3},\;
\;\;\;\int_0^{\infty}\frac{64t^5}{n^6+exp(t^6)}dt=\frac{32ln(n^6+1)}{3n^6}\hspace{.1in}.\;$$
So we arrive at the first assertion of the theorem.

 Now for the proof of the second part of the theorem, that is for the
case
 $k=o(n^{3/2})$  we study the asymptotic behavior of $ES_u(a,b)$  for different intervals
 $(-\infty,-1)$, $(-1,0)$, $(0,1)$ and $(1,\infty)$ separately.

 For   $1<x<\infty$,  using the change of variable
$x=1+\frac{t}{n}$, and by (1.4),(1.5), and  (2.1) we find that
\begin{eqnarray*}
\frac{1}{n}\int_{0}^{\infty}f_{n}(1+\frac{t}{n})dt=\frac{1}{\pi}\int_{0}^{\infty}
R_{1}(t)dt+o(1)
\end{eqnarray*}
where$\int_{0}^{\infty}R_{1}(t)dt=0.734874192$ .

 For
$-\infty<x<-1$  using   the change of variable $x=-1-\frac{t}{n}
$,and by  the fact  that $u^2/n^3=o(1)$, as $n\rightarrow
\infty$, we have  that
$$exp\{u^2G_{2,1}(-1-t/n)/n^3 \}=1+o(1)$$ Therefore the relations
(2.4)implies that  $\;exp\{ g_{_2}(-1-t/n)\}=1+o(1)$. Thus by the
relations (1.4),(1.5),(2.3),(2.4) we have that
\begin{eqnarray*}
\frac{1}{n}\int_{0}^{\infty}f_{n}(-1-\frac{t}{n})dt=\frac{1}{\pi}\int_{0}^{\infty}
R_{2}(t)dt+o(1)
\end{eqnarray*}
where $\int_{0}^{\infty}R_{2}(t)dt=1.095640061$.\\
For $0<x<1$  using the change of variable $x=1-\frac{t}{n+t}$, and
relations (1.4),(1.5),(2.5),(2.6), and by using the equality (2.7)
 we have  that
\begin{eqnarray*}
\frac{n}{(n+t)^2}\int_{0}^{\infty}f_{n}(-1-\frac{t}{n+t})dt=
\frac{1}{\pi}\int_{0}^{\infty}(R_{3}(t)-\frac{I[t>1]}{2t})dt
+\frac{1}{2\pi}\log{(2n+1)}+o(1)
\end{eqnarray*}
where
$\int_{0}^{\infty}(R_{3}(t)-\frac{I[t>1]}{2t})dt=-.28977126$.\\
For $-1<x<0$,  using the change of variable $x=-1+\frac{t}{n+t}$,
   (2.9), and  by the same reasoning  as  above, the case $-\infty<x<-1$, we
have  that $\;exp\{ g_{_2}(-1+\frac{t}{n+t})\}=1+o(1)$ . Thus by
using  the relations (2.9),(2.10), and by using the equality (2.7)
 we have  that
\begin{eqnarray*}
\frac{n}{(n+t)^2}\int_{0}^{\infty}f_{n}(-1+\frac{t}{n+t})dt=
\frac{1}{\pi}\int_{0}^{\infty}(R_{4}(t)-\frac{I[t>1]}{2t})dt
+\frac{1}{2\pi}\log{(2n+1)}+o(1)
\end{eqnarray*}
where$\int_{0}^{\infty}(R_{4}(t)-\frac{I[t>1]}{2t})dt=.3793914850$.
This complete the proof of the theorem.


\begin{thebibliography}{99}
\bibitem{bhs}  A  T  Bharucha-Ried  and  M  Sambandham,  {\em Random
  Polynomials}.  Academic Press, New York (1986).
\bibitem{crl}
 H  Cram\'{e}r  and  M  R  Leadbetter,
{\em      Stationary and Related Stochastic Processes},
  John  Wiley  New York (1967).

\bibitem{far1}  K  Farahmand,  {\em Real zeros of random algebraic
polynomials}, Proc. Amer. Math. Soc., 113(1991), 1077-1084.
\bibitem{far2} K  Farahmand,  (1998).  {\em Topics  in Random Polynomials},
Chapman \& Hall
\bibitem{ibm}  I  A  Ibragimov  and N  B   Maslova,
{\em On the expected number of real  zeros of random  polynomials,
coefficients with zero means}. Theory Probab. Appl.,  16(1971,
228-248.
\bibitem{kac1}  M Kac,   {\em On the average number of real  roots  of a
random algebraic equation}, Bull. Amer. Math. Soc., 49(1943),
314-320.
\bibitem{lsh1}  B F Logan  and  L A Shepp,  {\em  Real zeros of random
polynomials}, Proc. London Math. Soc., Ser. 3, 18(1968), 29-35.
\bibitem{lsh2}  B F Logan  and  L  A  Shepp,   {\em Real zeros of random
polynomials},  Proc. London Math. Soc., Ser. 3, 18(1968), 308-314.
\bibitem{rs3}  S  Rezakhah  and  S  Shemehsavar, {\em  On the
average  Number of Level Crossings of  Certain Gaussian Random
Polynomials}. Nonlinear. Anal.,  63,  (2005), pp. e555-567.
\bibitem{rs1}  S Rezakhah,  and A R Soltani,.  {\em On the Levy and Harmonizable random polynomials}.
 Georgian. Math. Journal, vol. 7, No. 2,(2000) pp.  379-386.
\bibitem{rs2}  S Rezakhah, A R  Soltani, {\em On the
Expected Number of Real Zeros of Certain Gaussian  Random
Polynomials}. Stoc. Anal, Appl, vol. 21, No. 1, pp. 223-234.
\bibitem{sam1}  M  Sambandham,  {\em On the real roots of the random algebraic
 polynomial},  Indian J. Pure Appl. Math., 7(1976), 1062-1070.
\bibitem{sam2}  M  Sambandham, {\em  On a random algebraic
 equation},  J. Indian  Math. Soc., 41(1977), 83-97.
\bibitem{wil}  J  E  Wilkins,  {\em  An asymptotic expansion
for the expected number of real zeros  of a random polynomial}.
Proc. Amer. Math. Soc. 103(1988), 1249-1258.
\end{thebibliography}
\end{document}